\documentclass[a4paper,reqno]{amsart}
\usepackage{amssymb}
\title{Generalized Lawson tori and Klein bottles}
\author{Alexei V. Penskoi}
\address{Department of Higher Geometry and Topology, 
Faculty of Mathematics and Mechanics, Moscow State University,
Leninskie Gory, GSP-1, 119991, Moscow, Russia\newline \emph{and}
\newline Independent University of Moscow, 
Bolshoy Vlasyevskiy pereulok~11, 119002, Moscow, Russia
\newline \emph{and}
\newline Faculty of Mathematics, Higher School of Economics,
Vavilova Str.~7, 117312, Moscow, Russia}
\email{penskoi@mccme.ru}
\date{}
\newtheorem{Proposition}{Proposition}
\newtheorem{Definition}{Definition}
\newtheorem{Theorem}{Theorem}
\DeclareMathOperator{\Area}{Area}
\DeclareMathOperator{\sn}{sn}
\DeclareMathOperator{\cn}{cn}
\DeclareMathOperator{\dn}{dn}
\DeclareMathOperator{\Ec}{Ec}
\DeclareMathOperator{\Es}{Es}
\DeclareMathOperator{\am}{am}
\begin{document}
\begin{abstract}
Using Takahashi theorem we propose an approach
to extend known families of minimal tori in
spheres. As an example, the well-known two-parametric
family of Lawson tau-surfaces including
tori and Klein bottles is extended to a three-parametric
family of tori and Klein bottles minimally immersed
in spheres. Extremal spectral
properties of the metrics on these surfaces are investigated.
These metrics include i) both metrics extremal for the first
non-trivial eigenvalue on the torus, i.e. the metric on the
Clifford torus and the metric on the equilateral torus
and ii) the metric maximal for the first non-trivial eigenvalue on
the Klein bottle.

\end{abstract}
\maketitle

\section{Introduction}

\subsection{The Lam\'e equation and the statement of the main Theorem}

The well-known Lam\'e equation is usually written as
\begin{equation}\label{Lame}
\frac{d^2\varphi}{dz^2}+(h-n(n+1)(k\sn z)^2)\varphi=0,
\end{equation}
where $k$ is the module of $\sn z,$
see e.g. the books~\cite{Arscott1964,BatemanErdelyi1955}.

In the case $n=1$ three
wonderful solutions of the Lam\'e equation~\eqref{Lame}
given by three Jacobi elliptic functions are known,
$$
\Ec_1^0(z)=\dn z,\quad\Ec_1^1(z)=\cn z,\quad\Es_1^1(z)=\sn z,
$$
where we use the Ec/Es notation for the solutions
used by Ince in the paper~\cite{Ince1940}.
These solutions correpond to
\begin{equation}\label{h-values}
h=k^2, \quad h=1,\quad h=1+k^2
\end{equation}
respectively.

The change of variable
\begin{equation}\label{sn-sin}
\sn z=\sin y\quad\Longleftrightarrow\quad y=\am z,
\end{equation}
where $\am z$ is Jacobi amplitude function,
see e.g. the book~\cite[Section~13.9]{BatemanErdelyi1955},
transforms the Lam\'e equation~\eqref{Lame}
into a trigonometric form of the Lam\'e
equation
\begin{equation}\label{Lame-trig}
[1-(k\sin y)^2]\frac{d^2\varphi}{dy^2}-%
k^2\sin y\cos y\frac{d\varphi}{dy}+%
[h-n(n+1)(k\sin y)^2]\varphi=0.
\end{equation}

This trigonometric form of the Lam\'e equation
is used in the book~\cite{Arscott1964}.
The change of variable $\sn z=\cos y$
leads to another trigonometric form
used in the book~\cite{BatemanErdelyi1955}.

Using standard properties of the Jacobi elliptic functions
and the change of variable~\eqref{sn-sin} one obtains
three solutions of the Lam\'e equation in the trigonometric
form~\eqref{Lame-trig},
\begin{equation}\label{Lame-solutions}
\Ec_1^0(y)=\sqrt{1-k^2\sin^2y},\quad%
\Ec_1^1(y)=\cos y,\quad\Es_1^1(y)=\sin y.
\end{equation}

Let us consider functions
\begin{gather}
\tilde\varphi_1(y)=\sqrt{\frac{b^2+c^2-a^2}{2(c^2-a^2)}}\sin y, 
\quad \tilde\varphi_2(y)=\sqrt{\frac{a^2+c^2-b^2}{2(c^2-b^2)}}\cos y,\label{phi12}\\
\tilde\varphi_3(y)=\sqrt{\frac{a^2+b^2-c^2}{2(b^2-c^2)}}\sqrt{1-\frac{b^2-a^2}{c^2-a^2}\sin^2y}.\label{phi3}
\end{gather}
These functions are rescaled three solutions~\eqref{Lame-solutions} 
of the Lam\'e equation in trigonometric 
form~\eqref{Lame-trig}
with $n=1$ and $k=\sqrt{\frac{b^2-a^2}{c^2-a^2}}.$

Let us denote by $K(\cdot)$ and $E(\cdot)$ the complete
elliptic integrals of the first and second kind
respectively defined as in the book~\cite{BatemanErdelyi1955}
by formulae
$$
K(k) = \int\limits_0^1\frac{d\alpha}{\sqrt{1-\alpha^2}\sqrt{1-k^2\alpha^2}},
\quad E(k)= \int\limits_0^1\frac{\sqrt{1-k^2\alpha^2}}{\sqrt{1-\alpha^2}}d\alpha.
$$

The main result of this paper is the following theorem.

\begin{Theorem}\label{maintheorem}
Let 
$F_{a,b,c}:\mathbb{R}^2\longrightarrow\mathbb{S}^5\subset\mathbb{R}^6$
be a three-parametric
doubly-periodic immersion of the plane to the 5-dimensional sphere of radius $1$
defined by the formula
\begin{gather}
F_{a,b,c}(x,y)=(\sin ax\,\tilde\varphi_1(y),\cos ax\,\tilde\varphi_1(y),\label{Fabc}\\
\sin bx\,\tilde\varphi_2(y),\cos bx\,\tilde\varphi_2(y),
\sin cx\,\tilde\varphi_3(y),\cos cx\,\tilde\varphi_3(y)),\notag
\end{gather}
where
\begin{itemize}
\item[a)] either $a,b,c$ are integers and $|c|>\sqrt{a^2+b^2},$
\item[b)] or $a,b$ are nonzero integers and $|c|=\sqrt{a^2+b^2}.$
\end{itemize}

Let $\mathcal{L}=\{(2\pi n,2\pi m)|n,m\in\mathbb{Z}\}$
and $\tilde{F}_{a,b,c}:\mathbb{R}^2/\mathcal{L}\longrightarrow\mathbb{S}^5\subset\mathbb{R}^6$
be the natural map induced by $F_{a,b,c}.$

Let $S(a,b,c)=\frac{4\pi}{\sqrt{c^2-a^2}}%
\left(2(c^2-a^2)E\left(\sqrt{\frac{b^2-a^2}{c^2-a^2}}\right)%
-(c^2-a^2-b^2)K\left(\sqrt{\frac{b^2-a^2}{c^2-a^2}}\right)\right).$

Then the following statements hold.

1) The image $T_{a,b,c}=F_{a,b,c}(\mathbb{R}^2)$
is a minimal compact surface in the 5-dimensional sphere
$(\mathbb{S}^5).$

2) The case b) corresponds to Lawson tau-surfaces 
$\tau_{a,b}\cong T_{a,b,\sqrt{a^2+b^2}}.$
Distinct Lawson tau-surfaces correpond to unordered pairs
$a,b\geqslant1$ such that $(a,b)=1.$
The surface $T_{a,b,\sqrt{a^2+b^2}}$ is a Lawson torus $\tau_{a,b}$
if $a$ and $b$ are odd and $T_{a,b,\sqrt{a^2+b^2}}$ is a Lawson Klein
bottle $\tau_{a,b}$ if either $a$ or $b$ is even, where
we assume $(a,b)=1.$

3) In the case b) the metric induced on
$\tau_{a,b}\cong T_{a,b,\sqrt{a^2+b^2}}$
is extremal for the functionals
$\Lambda_j(\mathbb{T}^2,g)$ if $\tau_{a,b}$
is a Lawson torus or $\Lambda_j(\mathbb{KL},g)$
if $\tau_{a,b}$ is a Lawson Klein bottle,
where
$j=2\left[\frac{\sqrt{a^2+b^2}}{2}\right]+a+b-1$
and $[\cdot]$ denotes the integer part.
The corresponding value of the functional is
$\Lambda_j(\tau_{a,b})=8\pi aE\left(\frac{\sqrt{a^2-b^2}}{a}\right).$

4) In the case a) for an integer $k\geqslant1$ one has
$T_{a,b,c}=T_{ka,kb,kc}.$ Moreover,
$T_{-a,b,c},$ $T_{a,-b,c}$ $T_{a,b,-c}$ and $T_{b,a,c}$
are isometric to $T_{a,b,c}.$
Hence, it is sufficient to consider non-negative integer $a,b,c$
satisfying conditions a) such that $(a,b,c)=1$ and assume that $(a,b,c)$
and $(b,a,c)$ are equivalent.

5) In the case a) depending on the parity of $a,$ $b$ and $c$ we have the following three 
subcases.
\begin{itemize}
\item[I)] If $a$ and $b$ have different parity and $c$ is even then
the surface $T_{a,b,c}$ is a Klein bottle and 
$\tilde{F}_{a,b,c}:\mathbb{R}^2/\mathcal{L}\longrightarrow T_{a,b,c}$
is a double covering. The area of $T_{a,b,c}$ is equal to $\frac{1}{2}S(a,b,c).$
\item[II)] If $a$ and $b$ are odd and $c$ is even then
the surface $T_{a,b,c}$ is a torus and 
$\tilde{F}_{a,b,c}:\mathbb{R}^2/\mathcal{L}\longrightarrow T_{a,b,c}$
is a double covering. The area of $T_{a,b,c}$ is equal to $\frac{1}{2}S(a,b,c).$
\item[III)] Otherwise, the surface $T_{a,b,c}$ is a torus and 
$\tilde{F}_{a,b,c}:\mathbb{R}^2/\mathcal{L}\longrightarrow T_{a,b,c}$
is a one-to-one map. The area of $T_{a,b,c}$ is equal to $S(a,b,c).$
\end{itemize}

6) In the case a) the metric induced on the torus or the Klein bottle
$T_{a,b,c}$ is extremal for
the functional $\Lambda_j(\mathbb{T}^2,g)$ or $\Lambda_j(\mathbb{KL},g)$
respectively, where
\begin{itemize}
\item[I)] if $a$ and $b$ have different parity and $c$ is even then
$j=a+b+c-3$ except the case of $T_{a,0,c}$
where $j=a+c-2,$
\item[II)] if $a$ and $b$ are odd and $c$ is even then $j=a+b+c-3,$
\item[III)] otherwise, $j=2(a+b+c)-3$ except the case of $T_{a,0,c}$
where $j=2(a+c)-2$ and the case of $T_{0,0,1}$
where $j=1.$
\end{itemize}

The corresponding value of this functional $\Lambda_j(T_{a,b,c})$
is $S(a,b,c)$ in the subcases I) and II) and 
$2S(a,b,c)$ in the subcase III).

\end{Theorem}

Remark that in the case a) $a$ and $b$ could be zero.
In particular, $T_{0,0,1}$ is the Clifford torus but
with the metric multiplied by $\frac{1}{2}.$ Hence,
it is not suprising that the metric on
$T_{0,0,1}$ is extremal for $\Lambda_1(\mathbb{T}^2,g).$

We will also explain in Section~\ref{T102-section}
that $T_{1,0,2}$ turns out to be isometric to 
Klein bottle $\tilde{\tau}_{3,1}$
carrying the metric extremal for the first eigenvalue
on the Klein bottle.

Let us now explain the statement of this theorem,
the idea of extending families of minimal tori and the formal proof.

\subsection{Minimal tori in spheres}\label{minimal-tori-in-spheres}

In his 1970 paper~\cite{Lawson1970} Lawson introduced several families
of minimal surfaces in $S^3$ including a family $\tau_{m,n}.$

\begin{Definition}\label{Lawson1970}
A Lawson tau-surface $\tau_{m,n}\looparrowright\mathbb{S}^3$
is defined as the image of the doubly-periodic immersion 
$\Psi_{m,n}:\mathbb{R}^2\looparrowright\mathbb{S}^3\subset\mathbb{R}^4$
given by the explicit formula
\begin{equation}\label{immersion}
\Psi_{m,n}(x,y)=%
(\cos(mx)\cos y, \sin(mx)\cos y, \cos(nx)\sin y,\sin(nx)\sin y).
\end{equation}
\end{Definition}

Here and later $\looparrowright$ denotes an immersion.

Lawson proved that for each unordered pair 
of positive integers $(m,n)$ with $(m,n)=1$ 
the surface $\tau_{m,n}$ is a distinct compact
minimal surface in $\mathbb{S}^3.$
Let us impose the condition $(m,n)=1.$
If both integers
$m$ and $n$ are odd then $\tau_{m,n}$ is a torus.
We call it a Lawson torus. If
one of integers $m$ and $n$ is even then $\tau_{m,n}$ is a
Klein bottle. We call it a Lawson Klein bottle.
Remark that  $m$ and $n$ cannot both
be even due to the condition $(m,n)=1$. 
The torus $\tau_{1,1}$ is the Clifford torus.

As explained in the statement 2) of Theorem~\ref{maintheorem},
the surfaces $T_{a,b,c}$ introduced in the Theorem~\ref{maintheorem}
are generalizations of the Lawson tau-surfaces.

Since Lawson paper~\cite{Lawson1970} several methods 
for constructing or describing minimal
tori in spheres were developed. An exhaustive review of all
such methods requires writing a book, hence we can mention here
only several ones.

Hsiang and Lawson developed in their 
paper~\cite{Hsiang-Lawson1971} a theory of reduction of a minimal
submanifold by a group action. This theory reduces the question
about construction of $\mathbb{S}^1$-invariant minimal tori
to the question about construction of closed geodesics,
which is much simpler. As
the simplest example of application of this approach 
one can consider a family of Otsuki tori $O_{\frac{p}{q}}$
minimally immersed in $\mathbb{S}^3.$ They were introduced by Otsuki 
in his paper~\cite{Otsuki1970} using another approach,
but a detailed treatment of the Otsuki tori
using Hsiang-Lawson approach could be found in the
paper~\cite{Penskoi2013}. Hsiang and Lawson briefly
mentioned Otsuki tori but they gave also new
examples of $\mathbb{S}^1$-invariant minimal tori in $\mathbb{S}^3.$
As another example one can mention the
paper~\cite{Ferus-Pedit1990} by Ferus and Pedit where
Hsiang-Lawson approach is applied in the case of
$\mathbb{S}^1$-invariant minimal tori in $\mathbb{S}^4.$

Unfortunately, this approach gives a description of
families of minimal tori but does not give explicit
formulae for these tori. For example, Otsuki
tori $O_{\frac{p}{q}}$ are in one-to-one correspondence
with rational numbers $\frac{p}{q}$
such that $\frac{1}{2}<\frac{p}{q}<\frac{\sqrt{2}}{2},$ $p,q>0,$ $(p,q)=1,$
but we do not know explicit formulae for Otsuki
tori since the reconstruction of the torus $O_{\frac{p}{q}}$
from a fraction $\frac{p}{q}$ requires solving a transcendental equation
and a system of ODEs.

Another approach is based on methods of integrable systems
and describes minimal tori in spheres through
algebraic geometry. This approach was developed by many authors
in different particular cases starting from the 
paper~\cite{Hitchin1990} by Hitchin dealing with the case of 
$\mathbb{S}^3$ and finishing by the paper~\cite{Burstall1995} 
by Burstall dealing
with the general case of $\mathbb{S}^n.$ In fact, the investigation
of minimal tori in spheres was a part of an extended study by many
authors of harmonic maps from tori into symmetric spaces.
Let us mention here e.g. the paper~\cite{Ferus-Pedit-Pinkall-Sterling1992}
by Ferus, Pedit, Pinkall and Sterling dealing with the case
of $\mathbb{S}^4.$
We refer the reader to the
recent paper~\cite{Carberry2012} by Carberry containing a review
of the current situation of this approach with an extended list
of references.

This method describes all minimal tori in $\mathbb{S}^n$
through data including an algebraic
curve of genus $\gamma,$ a divisor $\mathcal{D}=P_1+\dots+P_\gamma$
consisting of $\gamma$ points and some additional data satisfying
so called {\em periodicity conditions.} The good news is that
a minimal torus can be reconstructed from these algebro-geometric
spectral data through complicated but in principle
explicit formulae involving theta-functions
of genus $\gamma.$ The bad news is that there is no constructive
description of algebro-geometric spectral data satisfying the periodicity
conditions.

In the papers~\cite{Carberry2007,Carberry2012} Carberry studied intensively
these periodicity conditions but provided only existence results. We would like
to cite here one of them interesting for our goals.

\begin{Theorem}[Carberry, \cite{Carberry2007}]\label{carberrytheorem}
For each integer $n\geqslant0$, there are countably many real $n$-dimensional
families of minimal immersions from rectangular tori to $\mathbb{S}^3.$ 
Each family consists of maps from a fixed torus.
\end{Theorem}

There exist also several other approaches generating particular
examples of minimal surfaces in spheres. We would like to
mention here an approach by Mironov for constructing Hamiltonian-minimal
Lagrangian embeddings in $\mathbb{C}^N$ based on intersections 
of real quadrics of special type, see the paper~\cite{Mironov2004}.
As a by-product of this construction one obtains minimal surfaces
in spheres, see the paper~\cite{Karpukhin2013q1} by Karpuhin
for more details. This approach is interesting for us since it provides
a two-parametric family $M_{m,k}$
of tori minimally immersed in $\mathbb{S}^5.$ It could be easily verified
that this family is a subfamily of our family from Theorem~\ref{maintheorem},
$M_{m,k}\cong T_{m,k,m+k}.$ This family was described in conformal coordinates
in the paper~\cite{Haskins2004} by Haskins and in the paper~\cite{Joyce2002} by
Joyce, but it seems that this family first appeared in a parametrization
similar to $T_{m,k,m+k}$ in the paper~\cite{Mironov2010} by Mironov. One can find
a detailed study of $M_{m,k}$
in the paper~\cite{Karpukhin2013q1} by Karpukhin .

Thus, minimal tori in spheres are described but in implicit way and only
several explicitly parametrized examples are known. However, recent
progress in study of extremal metrics brings minimal surfaces
in spheres back to our attention.

\subsection{Extremal metrics and minimal surfaces in spheres}

Let $M$ be a closed surface and $g$ be
a Riemannian metric on $M.$
Let us consider the associated Laplace-Beltrami
operator $\Delta:C^\infty(M)\longrightarrow C^\infty(M),$
$$
\Delta f=-\frac{1}{\sqrt{|g|}}\frac{\partial}{\partial x^i}%
\left(\sqrt{|g|}g^{ij}\frac{\partial f}{\partial x^j}\right),
$$
and its eigenvalues
$$
0=\lambda_0(M,g)<\lambda_1(M,g)\leqslant%
\lambda_2(M,g)\leqslant\lambda_3(M,g)\leqslant\dots
$$
Since the eigenvalues possess the following rescaling property,
$$
\forall t>0\quad\lambda_i(M,tg)=\frac{\lambda_i(M,g)}{t},
$$
it is natural to consider ``normalized'' eigenvalues
$$
\Lambda_i(M,g)=\lambda_i(M,g)\Area(M,g)
$$
invariant under 
the rescaling transformation $g\mapsto tg.$

Let us fix the surface $M$ and consider
$\Lambda_i(M,g)$ as a functional
$g\mapsto\Lambda_i(M,g)$ on the space
of all Riemannian metrics on $M.$

It turns out that the question
about the supremum $\sup\Lambda_i(M,g)$
of the functional $\Lambda_i(M,g)$ over the space
of Riemannian metrics $g$ on a fixed surface
$M$ is very difficult and only few results are known.

It is known that this supremum is finite
since functionals $\Lambda_i(M,g)$ are bounded from above.
It was proven in the paper~\cite{Yang-Yau1980} by Yang and Yau that
for an orientable surface $M$ of genus $\gamma$ the following
inequality holds,
$$
\Lambda_1(M,g)\leqslant 8\pi(\gamma+1).
$$
Korevaar proved in the paper~\cite{Korevaar1993}
that there exists a constant $C$ such that for
any $i>0$ and any compact surface $M$ of genus $\gamma$
the functional $\Lambda_i(M,g)$ is bounded,
$$
\Lambda_i(M,g)\leqslant C(\gamma+1)i.
$$

\begin{Definition} A metric $g_0$ on a fixed surface $M$
is called maximal for the functional $\Lambda_i(M,g)$
if
$$
\sup\Lambda_i(M,g)=\Lambda_i(M,g_0),
$$
where the supremum is taken
over the space
of Riemannian metrics $g$ on the fixed surface
$M.$
\end{Definition}

Only few maximal metrics are known at this moment.
The maximal metric for $\Lambda_1(\mathbb{S}^2,g)$
is the standard metric on the sphere (Hersch,~\cite{Hersch1970}),
the maximal metric for $\Lambda_1(\mathbb{RP}^2,g)$
is the standard metric on the
projective plane (Li and Yau,~\cite{Li-Yau1982}),
the maximal metric for $\Lambda_1(\mathbb{T}^2,g)$
is the metric on equilateral torus
(Nadirashvili,~\cite{Nadirashvili1996}). 
The last known (and quite surprising) maximal metric 
is the maximal metric for the 
first eigenvalue $\Lambda_1(\mathbb{K}\mbox{l},g)$
on the Klein bottle. As it was proved in
El Soufi, Giacomini and Jazar paper~\cite{ElSoufi-Giacomini-Jazar2006}
using results of Jakobson, Nadirashvili and 
Polterovich paper~\cite{Jakobson-Nadirashvili-Polterovich2006},
this is the metric on the bipolar
Lawson surface $\tilde{\tau}_{3,1}.$

We know also an example where $\sup\Lambda_i(M,g)$ is known,
but however there is no (smooth) maximal metric. It was proved
by Nadirashvili in the paper~\cite{Nadirashvili2002}
that $\sup\Lambda_2(\mathbb{S}^2,g)=16\pi$
and the maximum is reached on a singular metric which can be obtained
as the metric on the union of two touching
spheres of equal radius with
canonical metric.

If one would like to find a maximum of a function of several variables,
then one usually starts by finding extrema of this function. The same
idea is also reasonable for the functionals $\sup\Lambda_i(M,g).$
However, one should be careful here.
The functional $\Lambda_i(M,g)$ depends continuously
on the metric $g,$ but this functional
is not differentiable. However, for analytic
deformations $g_t$ the left and right derivatives
of the functional $\Lambda_i(M,g_t)$ with respect to $t$ exist,
see the papers by Berger~\cite{Berger1973}, Bando and Urakawa~\cite{Bando-Urakawa1983},
El~Soufi and Ilias~\cite{ElSoufi-Ilias2008}.
This led  to the following definition,
see the paper~\cite{Nadirashvili1996}
by Nadirashvili and the papers~\cite{ElSoufi-Ilias2000,ElSoufi-Ilias2008}
by El Soufi and Ilias.

\begin{Definition}
A Riemannian metric $g_0$ on  a closed surface
$M$ is called extremal metric for the
functional $\Lambda_i(M,g)$ if for any analytic deformation
$g_t$ the following inequality holds,
$$
\frac{d}{dt}\Lambda_i(M,g_t)%
\left.\vphantom{\raisebox{-0.5em}{.}}\right|_{t=0+}\cdot%
\frac{d}{dt}\Lambda_i(M,g_t)%
\left.\vphantom{\raisebox{-0.5em}{.}}\right|_{t=0-}\leqslant0.
$$
\end{Definition}

Investigation of extremal metrics turned out to be useful.
For example, Jakobson, Nadirashvili and 
Polterovich proved in the 
paper~\cite{Jakobson-Nadirashvili-Polterovich2006}
that the mentioned above metric on the Klein bottle realized
as the bipolar Lawson surface $\tilde{\tau}_{3,1}$
is extremal for $\Lambda_1(\mathbb{KL},g)$
and using this result El Soufi, Giacomini and 
Jazar proved in the paper~\cite{ElSoufi-Giacomini-Jazar2006}
the above mentioned result
that this metric is the unique extremal metric and the maximal
one.

As one can expect, we know more about extremal metrics then about
maximal metrics. El~Soufi
and Ilias proved in the paper~\cite{ElSoufi-Ilias2000}
that the only extremal
metric for $\Lambda_1(\mathbb{T}^2,g)$ different from the maximal one
is the metric on the Clifford torus.

Let us remark that the metrics on the surfaces $T_{a,b,c}$ from Theorem~\ref{maintheorem}
includes both metrics extremal for the first
eigenvalue on the torus, i.e. the metric on the
Clifford torus $\tau_{1,1}\cong T_{1,1,\sqrt{2}}$
and the metric on the equilateral torus $M_{1,1}\cong T_{1,1,2}.$
As we will show in section~\ref{T102-section}, these metrics also
include the metric on the Klein bottle $T_{1,0,2}\cong\tilde{\tau}_{3,1}$
extremal for the $\Lambda_1(\mathbb{KL},g).$ Hence, $T_{a,b,c}$
includes surfaces carrying all extremal metrics for the first
eigenvalue on the torus and Klein bottle.

Extremality of several families of metrics on the torus and Klein bottles
was investigated recently. 
\begin{itemize}
\item Lapointe studied metrics
on bipolar Lawson surfaces $\tilde{\tau}_{m,k}\looparrowright\mathbb{S}^4$
in his 2008 paper~\cite{Lapointe2008}.
\item The author studied metrics on Lawson surfaces $\tau_{m,k}\looparrowright\mathbb{S}^3$
and metrics on Otsuki tori $O_{\frac{p}{q}}\looparrowright\mathbb{S}^3$
in his 2012 paper~\cite{Penskoi2012}
and 2013 paper~\cite{Penskoi2013} respectively.
\item Karpukhin studied metrics
on bipolar Otsuki tori $\tilde{O}_{\frac{p}{q}}\looparrowright\mathbb{S}^4$
and on the family of tori $M_{m,k}\looparrowright\mathbb{S}^5$
in his 2013 paper~\cite{Karpukhin2013}
and the paper~\cite{Karpukhin2013q1} respectively.
\item Karpukhin also proved in the paper~\cite{Karpukhin2013q2}
that all metrics mentioned in this list are not maximal
except metrics on
$M_{1,1}$ (the equivalateral torus) and $\tilde{\tau}_{3,1}.$
\end{itemize}

The significant progress in study of extremal metrics in
the papers~\cite{Karpukhin2013,Karpukhin2013q1,Karpukhin2013q2,Penskoi2012,Penskoi2013}
became possible due to El Soufi-Ilias theorem establishing
relation between extremal metrics and minimal surfaces in spheres.

Let $M$ be a two-dimensional minimally immersed
submanifold of the standard
sphere $\mathbb{S}^n\subset\mathbb{R}^{n+1}$
of radius~$1.$ Let $\Delta$ be the Laplace-Beltrami operator
on $M$ equipped with the induced metric.

Let us introduce the Weyl eigenvalues counting function
$$
N(\lambda)=\#\{\lambda_i|\lambda_i<\lambda\}.
$$
Remember that we count the
eigenvalues starting from $\lambda_0=0.$

\begin{Theorem}[El Soufi, Ilias, \cite{ElSoufi-Ilias2008}]\label{ElSoufi-Ilias}
The metric induced on $M$ by the immersion
$M\looparrowright\mathbb{S}^n$ is an extremal metric
for the functional $\Lambda_{N(2)}(M,g).$
\end{Theorem}

Thus, investigation of (smooth) extremal metrics
on surfaces could be done in the following way:
\begin{itemize}
\item find a minimal surface in a sphere,
\item find $N(2),$
\item then the metric induced on the minimal surface is extremal
for the functional $\Lambda_{N(2)}.$
\end{itemize}
However, it is not easy to follow this approach.
As we discussed in Section~\ref{minimal-tori-in-spheres},
even the descriptions of minimal tori in spheres
are quite complicated and implicit. Moreover, it
turns out that there is no known general way to find
$N(2)$ and in each example one should
invent an ad hoc argument.

All mentioned above successful examples
of application of this approach in 
papers~\cite{Karpukhin2013,Karpukhin2013q1,Karpukhin2013q2,Penskoi2012,Penskoi2013}
share the following features:
\begin{itemize}
\item these surfaces were already known to be minimal in spheres,
\item their metrics are metrics of revolution,
\item either their parametrisation is explicit (Lawson surfaces
and Lawson bipolar surfaces) or the structure of zeroes
of immersion functions is simple (Otsuki tori and Otsuki
bipolar tori).
\end{itemize}

At this moment there is no hope to investigate
all extremal metrics on tori since this requires
at least a constructive description of minimal tori
in spheres and it seems that the existing 
implicit description in terms of algebro-geometric
data could not be improved. Thus we concentrate now
our efforts on investigating particular examples
of extremal metrics. This leads us
to the problem of finding new explicit
examples of minimal tori in spheres.

\subsection{Constructing explicit examples
of minimal tori via Takahashi theorem}

Let us recall the well-known result about
description of the minimal surfaces
in $\mathbb{R}^n$ in terms of harmonic
functions.

\begin{Proposition}
A submanifold $M\looparrowright\mathbb{R}^n$
is minimal if and only if the restictions
$x^1|_M,\dots x^n|_M$ to $M$ of the coordinate
functions in $R^n$ are harmonic with respect to
the Laplace-Beltrami operator $\Delta^M$ on $M$ equipped
with the induced metric,
$$
\Delta^M x^i|_M=0.
$$
\end{Proposition}

One can rewrite this Proposition in terms
of isometric immersions.

\begin{Proposition}
An isometric immersion $f:M\looparrowright\mathbb{R}^n$
is minimal if and only if the components
of the immersion $f=(f^1,\dots,f^n)$
are harmonic with respect to the Laplace-Beltrami operator 
$\Delta$ on $M,$
$$
\Delta f^i=0.
$$
\end{Proposition}

If an isometric immersion by harmonic functions
(i.e. eigenfunction of $\Delta$ with eigenvalue 0)
is minimal in $\mathbb{R}^n,$ what can we say about
isometric immersions by eigenfunctions of $\Delta$
with a common eigenvalue $\lambda?$ The answer is
given by the Takahashi theorem.

\begin{Theorem}[Takahashi,~\cite{Takahashi1966}]\label{takahashi}
An isometric immersion $f:M\looparrowright\mathbb{R}^{n+1},$
where $f=(f^1,\dots,f^{n+1}),$
defined by eigenfunctions
$f^i$ of the Laplace-Beltrami operator $\Delta$
with a common eigenvalue $\lambda,$
$$
\Delta f^i=\lambda f^i,
$$
possesses the following properties,
\begin{itemize}
\item the image $f(M)$ lies on the sphere $\mathbb{S}^n_R$
of radius $R$ with the center at the origin
such that 
\begin{equation}\label{lambda-R}
\lambda=\frac{\dim M}{R^2},
\end{equation}
\item the immersion $f:M\looparrowright\mathbb{S}^n_R$
is minimal.
\end{itemize}

If $f:M\looparrowright\mathbb{S}^{n}_R,$
where $f=(f^1,\dots,f^{n+1}),$ is a minimal
isometric immersion of a manifold $M$
into the sphere $\mathbb{S}^{n}_R$
of radius $R,$ then $f^i$
are eigenfunctions of the Laplace-Beltrami
operator $\Delta,$
$$
\Delta f^i=\lambda f^i,
$$
with the same eigenvalue $\lambda$
such that $\lambda=\frac{\dim M}{R^2}.$
\end{Theorem}

Takahashi theorem describes minimal immersions in terms
of eigenfunctions of the Laplace-Beltrami operator.
This is a system of PDEs equivalent to the standard
system of PDEs describing minimal immersions
in terms of mean curvature normal vector.
This description is more natural from the point
of view of spectral geometry.

Surprisingly, this approach did not generate
much interest till very recently. To the best
of the author's knowledge, the only known
successful application of this approach
to construction of minimal surfaces in spheres
was in the above mentioned
paper~\cite{Jakobson-Nadirashvili-Polterovich2006}.
In this paper using Takahashi theorem and
properties of eigenfunctions
Jakobson, Nadirashvili and Polterovich
constructed a minimal isometric immersion
of the Klein bottle to $\mathbb{S}^4$
such that the corresponding metric is
extremal for the first eigenvalue. This
metric turned out to be Lawson bipolar surface $\tilde{\tau}_{3,1}.$

However, solving a PDE system is a difficult task.
Is it possible to find at least some new particular examples
of minimal surfaces in spheres using the Takahashi theorem? 
In the present paper we propose an approach
based on the Takahashi theorem leading us to
an extension of known families of minimal tori in spheres
using solving systems of algebraic equations. On this way
we obtain Theorem~\ref{maintheorem}.

\section{Extension of families of minimal tori in spheres using
Takahashi theorem}\label{ansatz}

Let us start with a very naive idea. Let $x,y$
be coordinates in the plane $\mathbb{R}^2.$
Let us choose randomly a second order 
elliptic differential operator $L$ on the
plane invariant with respect to the translations
$(x,y)\mapsto (x+2\pi,y),$ $(x,y)\mapsto (x,y+2\pi).$
Consider spectral problem for $L$ with periodic
boundary conditions
$$
\left\{\begin{array}{l}
L\psi=\lambda\psi,\\
\psi(x,y)=\psi(x+2\pi,y)=\psi(x,y+2\pi).
\end{array}\right.
$$
Let $0=\lambda_0<\lambda_1\leqslant\lambda_2\leqslant\dots$
be the spectrum of this spectral problem. Let us
choose an eigenvalue $\lambda_j$ and
several linearly independent eigenfunctions
$f_1,\dots,f_{n+1},$ corresponding to the eigenvalue $\lambda_j.$
Consider the map 
$$
f:\mathbb{R}^2/\mathcal{L}\longrightarrow\mathbb{R}^{n+1},
$$
where $f=(f_1,\dots,f_{n+1}).$ 

Let $g$ be the pullback $f^*g_0$ of the standard euclidean metric $g_0$
in $\mathbb{R}^{n+1}$ to the torus $\mathbb{R}^2/\mathcal{L}$ by the map $f.$

It follows from Takahashi Theorem~\ref{takahashi} that if
\begin{itemize}
\item[(a)] $g$ is a Riemannian metric (i.e. positive definite) and
\item[(b)] the Laplace-Beltrami operator for $g$ coincides with the initial
differential operator $L,$
\end{itemize}
then $f$ is a minimal isometric immersion of the torus 
$(\mathbb{R}^2/\mathcal{L},g)$ to a sphere $\mathbb{S}^n_R$
of radius $R$ such that $\lambda_j=\frac{2}{R^2}.$

This idea is naive since there is practically no chance that
starting from a randomly chosen operator $L$ one can satisfy
conditions (a) and (b). Hence one should find a way to
start with a smart choice of the initial operator $L.$

Then it is time to rememeber Carberry Theorem~\ref{carberrytheorem}.
This theorem is for the case of $\mathbb{S}^3$ but we can hope that in the
general case minimal tori in spheres also like to exist in families.
That's why the key idea is the following: \emph{start with the operator
$L$ already known to be the Laplace-Beltrami operator on a
minimal torus in a sphere.}

Let us now remember that known examples of extremal metrics
are metrics on tori of revolution or on quotients of
tori on revolution, e.g. Lawson tori~\cite{Penskoi2012},
Otsuki tori~\cite{Penskoi2013}, bipolar Lawson 
tori~\cite{Lapointe2008}, bipolar Otsuki tori~\cite{Karpukhin2013}
and the family $M_{m,k}$ considered in~\cite{Karpukhin2013q1}. Let us then
restrict our attention to tori of revolution.

Since by rescaling we can always restrict our attention to the
case of the sphere of radius 1 corresponding by formula~\eqref{lambda-R}
to $\lambda=2$, we are interested in
the equation 
\begin{equation}\label{L-spectral}
L\psi=2\psi
\end{equation}
 and its solutions
of the form 
$$
\psi(x,y)=\varphi(y)\sin mx\quad\mbox{or}\quad 
\psi(x,y)=\varphi(y)\cos mx.
$$
After separation of variables in PDE~\eqref{L-spectral}
one obtains a linear second order ODE 
\begin{equation}\label{spectral-reduced}
\varphi''(y)+a(m,y)\varphi'(y)+b(m,y)\varphi(y)=0
\end{equation}
with periodic boundary conditions
\begin{equation}\label{spectral-reduced-boundary}
\varphi(y+2\pi)=\varphi(y).
\end{equation}
Thus, we propose the following method for constructing
minimal tori of revolution in spheres.

\begin{itemize}
\item Consider an equation
\begin{equation}\label{spectral-reduced-ours}
\varphi''(y)+A(\nu,y)\varphi'(y)+B(\nu,\mu,y)\varphi(y)=0
\end{equation}
with a spectral parameter $\mu$
and (possibly) an additional parameter $\nu$
such that this equation is 
already known to appear after a separation of variables
in the spectral problem~\eqref{L-spectral}
for the Laplace-Beltrami operator on a
minimal torus in the unitary sphere.
\item Take solutions $\varphi_1(y),\dots,\varphi_l(y)$
of equation~\eqref{spectral-reduced-ours}
with periodic boundary
conditions~\eqref{spectral-reduced-boundary},
corresponding to $\mu=\mu_1,\dots,\mu=\mu_l,$ respectively.
These solutions depend on the parameter $\nu.$
\item Consider the map $f:\mathbb{R}^2/\mathcal{L}\longrightarrow\mathbb{R}^{2l}$
given by the formula
$$
f(x,y)=(c_1\varphi_1(y)\cos m_1x,c_1\varphi_1(y)\sin m_1x,\dots,
c_l\varphi_l(y)\cos m_lx,c_l\varphi_l(y)\sin m_lx),
$$
where $c_1,\dots,c_l$ are constants.
\item Let $g=f^*g_0,$ where $g_0$ is the standard euclidean
metric on $\mathbb{R}^{2l}.$
\item Find the Laplace-Beltrami operator $L$ for the metric $g.$
Remark that $L$ depends on $\nu,c_1,\dots,c_l,m_1,\dots,m_l.$
\item Separate variables in the equation $L\psi(x,y)=2\psi(x,y)$
and obtain the equation
\begin{gather}
\varphi''(y)+a_T(\nu,c_1,\dots,c_l,m_1,\dots,m_l;m,y)\varphi'(y)+\notag\\
+b_T(\nu,c_1,\dots,c_l,m_1,\dots,m_l;m,y)\varphi(y)=0.\label{spectral-reduced-new}
\end{gather}
\item Write down a system of \emph{algebraic} (i.e. not differential)
equations
\begin{equation}\label{algebraic}
\left\{\begin{array}{l}
g_1(\nu,c_1,\dots,c_l,m_1,\dots,m_l)=0,\\
\ldots\\
g_N(\nu,c_1,\dots,c_l,m_1,\dots,m_l)=0,\\
\end{array}\right.
\end{equation}
equivalent to the condition of coincidence of initial equation~\eqref{spectral-reduced-ours}
with equation~\eqref{spectral-reduced-new}, i.e.
\begin{equation}\label{coincidence}
\left\{\begin{array}{l}
A(\nu,y)\equiv a_T(\nu,c_1,\dots,c_l,m_1,\dots,m_l;m_1,y),\\
B(\nu,\mu_1,y)\equiv b_T(\nu,c_1,\dots,c_l,m_1,\dots,m_l;m_1,y),\\
\dots\\
A(\nu,y)\equiv a_T(\nu,c_1,\dots,c_l,m_1,\dots,m_l;m_l,y),\\
B(\nu,\mu_l,y)\equiv b_T(\nu,c_1,\dots,c_l,m_1,\dots,m_l;m_l,y),\\
\end{array}\right.
\end{equation}
and add the condition
that $g$ is positive definite.
The sign ``$\equiv$'' in system~\eqref{coincidence}
means ``equals identically with respect to $y.$''
\item If system of equations~\eqref{algebraic} and the condition
that $g$ is positive definite have a solution
$(\nu,c_1,\dots,c_l,m_1,\dots,m_l),$ then by Takahashi Theorem~\ref{takahashi}
the map $f$ with these values of parameters
is an isometric immersion and the image $f(\mathbb{R}^2)$
is minimal in the unitary sphere $\mathbb{S}^{2l-1}.$
\end{itemize}

Let us consider Lawson tau-surfaces
and equation~\eqref{spectral-reduced}
appearing after separation of variables in the spectral
problem for the Laplace-Beltrami operator.
As we know from the paper~\cite{Penskoi2012},
this equation is the Lam\'e
equation in trigonometric form~\eqref{Lame-trig}.
Let us then in order to give an example of the proposed method
apply the described above algorithm to the Lam\'e
equation in trigonometric form~\eqref{Lame-trig}.
This gives us the family of surfaces from Theorem~\ref{maintheorem}.

In fact, we should remark that in this example we modify a little bit
the proposed approach since the image of $f$ could be not only
a torus but also a quotient of a torus, e.g. a Klein bottle.

\section{Proof of Theorem~\ref{maintheorem}}

Let us apply the algorithm from Section~\ref{ansatz}
taking the Lam\'e equation in trigonometric 
form~\eqref{Lame-trig} as equation~\eqref{spectral-reduced-ours}
and its classical solutions~\eqref{Lame-solutions}
as solutions, i.e. $l=3$ and
$$
\varphi_1(y)=\sin y, \quad \varphi_2(y)=\cos y, \quad
\varphi_3(y)=\sqrt{1-k^2\sin^2y},
$$
where $k$ is the module of the Lam\'e equation.
Here $k$ plays the role of the additional
parameter $\nu$ in equation~\eqref{spectral-reduced-ours}
that we should also find
from system~\eqref{algebraic}. The parameter $h$
in the Lam\'e equaion plays the role of the spectral
parameter $\mu$ in equation~\eqref{spectral-reduced-ours}
and as the values $\mu_i$ we take values of $h$
from formulae~\eqref{h-values}.

Then we have
\begin{gather*}
f(x,y)=(c_1\sin y\,\cos m_1x,c_1\sin y\,\sin m_1x,
c_2\cos y\,m_2x,c_2\cos y\,\sin m_2x,\\
c_3\sqrt{1-k^2\sin^2y}\,\cos m_3x,c_3\sqrt{1-k^2\sin^2y}\,\sin m_3x).
\end{gather*}
If $c_1,c_2,c_3,m_1,m_2,m_3,k$ satisfy system~\eqref{algebraic}
then by Takahashi theorem~\ref{takahashi} the image of $f$
is on the unitary sphere. This condition is equivalent to the equation
$$
-2+c_1^2+c_2^2+2c_3^2-c_3^2k^2+(-c_1^2+c_2^2+c_3^2k^2)\cos 2y\equiv0.
$$
It follows that 
\begin{equation}\label{c1c2}
c_1^2=1-c_3^2+c_3^2k^2,\quad c_2^2=1-c_3^2.
\end{equation}
A straightforward calculation shows that the metric $g=f^*g_0$
is given by the formula
\begin{gather*}
g=\frac{1}{2}[m_2^2(1-c_3^2)+m_3^2c_3^2(2-k^2)+
m_1^2(1+c_3^2(k^2-1))+\\
+(m_2^2(1-c_3^2)+m_3^2c_3^2k^2+ 
m_1^2(c_3^2-c_3^2k^2-1))\cos 2y]dx^2+\\
+\frac{k^2-2-2c_3^2(k^2-1)-k^2\cos 2y}{2k^2\sin^2 y-2}dy^2,
\end{gather*}
where we already applied formulae~\eqref{c1c2} to eliminate $c_1$ and $c_2.$

A straightforward calculation shows that the first equation in system~\eqref{coincidence}
implies the equation
\begin{gather*}
((k^2-1)(m_2^2(c_3^2-1)^2-m_3^2c_3^4k^2)+
m_1^2(1+c_3^2(k^2-1))^2)(2-k^2+k^2\cos 2y)\equiv0.
\end{gather*}
It follows that
\begin{equation}\label{m3}
m_3^2=\frac{m_2^2(c_3^2-1)^2(k^2-1)+m_1^2(1+c_3^2(k^2-1))^2}{c_3^4k^2(k^2-1)}.
\end{equation}
We use this formula to eliminate $m_3$ from our system of equations.

In the same way we investigate the second equation in system~\eqref{coincidence}.
We do not write down this equation explicitly
since it is a huge expression. It
turns out that this equation implies that
\begin{equation}\label{k}
k^2=\frac{(m_1^2-m_2^2)(2c_3^2-1)}{2m_1^2c_3^2+m_2^2(1-2c_3^2)}.
\end{equation}
Hence, we have three-parametric family of solutions parametrized
by $m_1,$ $m_2$ and $c_3.$ However, it is more convenient
to parametrize by $m_1,$ $m_2$ and $m_3.$ One can eliminate $k$
from equation~\eqref{m3} using equation~\eqref{k}
and obtain the equation
\begin{equation}\label{m3new}
m_3^2=\frac{m_1^2+m_2^2-2m_2^2c_3^2}{1-2c_3^2}.
\end{equation}
Then one can find $c_3$ from equation~\eqref{m3new}
and obtain the formula
\begin{equation}\label{c3}
c_3^2=\frac{m_1^2+m_2^2-m_3^2}{2(m_2^2-m_3^2)}.
\end{equation}
Now we can substitute formula~\eqref{c3} in equations~\eqref{k},
\eqref{c1c2} and obtain the formulae
\begin{gather*}
c_1^2=\frac{-m_1^2+m_2^2+m_3^2}{2(-m_1^2+m_3^2)},\quad
c_2^2=\frac{m_1^2-m_2^2+m_3^2}{2(-m_2^2+m_3^2)},\\
c_3^2=\frac{m_1^2+m_2^2-m_3^2}{2(m_2^2-m_3^2)},\quad
k^2=\frac{m_1^2-m_2^2}{m_1^2-m_3^2}.
\end{gather*}
Let us now define $\tilde\varphi_i=m_i\varphi_i$
and rename $a=m_1,$ $b=m_2,$ $c=m_3$ for the sake
of simplicity. Extracting square roots one obtains
formulae~\eqref{phi12}, \eqref{phi3} and \eqref{Fabc}.

Now we should satisfy periodicity conditions and also
choose such $a,b,c$ that $F_{abc}$ is a real map. This
is equivalent to conditions
\begin{itemize}
\item[] $a$ is integer and $\frac{b^2+c^2-a^2}{2(c^2-a^2)}>0$ or
\item[] $a$ is arbitrary and $\frac{b^2+c^2-a^2}{2(c^2-a^2)}=0,$
\end{itemize}
and
\begin{itemize}
\item[] $b$ is integer and $\frac{a^2+c^2-b^2}{2(c^2-b^2)}>0$ or
\item[] $b$ is arbitrary and $\frac{a^2+c^2-b^2}{2(c^2-b^2)}=0,$
\end{itemize}
and
\begin{itemize}
\item[] $c$ is integer and $\frac{a^2+b^2-c^2}{2(b^2-c^2)}>0$ or
\item[] $c$ is arbitrary and $\frac{a^2+b^2-c^2}{2(b^2-c^2)}=0.$
\end{itemize}
We should also add the condition that the metric $g$ is positive definite.

The solution of this system of conditions is exactly written in phrases
a) and b) in the statement of Theorem~\ref{maintheorem}.
Direct check shows that in fact, the values of $a,$ $b$ and $c$
satisfying these conditions give
minimal immersion of $\mathbb{R}^2/\mathcal{L}$ to the sphere $\mathbb{S}^5$
and the spectral problem after a separation of variables transforms
into the Lam\'e equation.
This finishes the proof of the statement 1) of Theorem~\ref{maintheorem}.

In the case b) our surfaces are Lawson tau-surfaces.
One obtains this by direct calculation. If $c=\sqrt{a^2+b^2}$ 
then $\varphi_1(y)=\sin y,$
$\varphi_2(y)=\cos(y)$ and $\varphi_3(y)=0.$ It is sufficient
now to remark that i) four first entries of the vector $F_{a,b,\sqrt{a^2+b^2}}$
coincide with $\Psi_{b,a}$ from~\eqref{immersion}
and the remaining two entries are zeroes and ii) $\tau_{a,b}\cong\tau_{b,a}.$
Then the statements 2) and 3) of Theorem~\ref{maintheorem}
follows from the papers~\cite{Lawson1970,Penskoi2012}.

This finishes the investigation of the case b) and in the following we
consider only the case a).

The images $T_{a,b,c}=F_{a,b,c}(\mathbb{R}^2)$ could be isometric for
distinct triples $(a,b,c).$ It is clear that for any integer $k>0$
one has $T_{a,b,c}=T_{ka,kb,kc}$. One can also remark that 
$T_{a,b,c}$ and $T_{-a,b,c}$ are isometric since $T_{-a,b,c}$
is the image of $T_{a,b,c}$ under the reflection of the ambient space
$\mathbb{R}^6$ with respect to the hyperplane $x^1=0.$
Similar statements are true for $T_{a,-b,c}$ and
$T_{a,b,-c}.$

Let us denote by $R$ the isometry
$$
R(x^1,x^2,x^3,x^4,x^5,x^6)=(x^3,x^4,-x^1,-x^2,x^5,x^6)
$$
of the ambient space $\mathbb{R}^6.$
Then we have the
identity $R\circ F_{b,a,c}(x,y+\frac{\pi}{2})=F_{a,b,c}(x,y).$
It follows that $T_{b,a,c}$ is isometric to $T_{a,b,c}.$
This implies the statement 4) of Theorem~\ref{maintheorem}.

Next, let us remark that $\tilde{F}_{a,b,c}$
is not necessarily a one-to-one map of
$\mathbb{R}^2/\mathcal{L}$ on $T_{a,b,c}$
because there could exist a non-trivial
(i.e. different from shifts by $2\pi$)
map $\Phi:\mathbb{R}^2\longrightarrow\mathbb{R}^2$
such that 
\begin{equation}\label{Phi}
F_{a,b,c}\circ\Phi=F_{a,b,c}.
\end{equation}
Let $(x_2,y_2)=\Phi(x_1,y_1),$ then equation~\eqref{Phi}
is equivalent to the system of equations
\begin{gather}
\sin ax_1\,\sin y_1=\sin ax_2\,\sin y_2,\label{ssn}\\
\cos ax_1\,\sin y_1=\cos ax_2\,\sin y_2,\label{csn}\\
\sin bx_1\,\cos y_1=\sin bx_2\,\cos y_2,\label{scn}\\
\cos bx_1\,\cos y_1=\cos bx_2\,\cos y_2,\label{ccn}\\
\sin cx_1\,\sqrt{1-k^2\sin^2y_1}=\sin cx_2\,\sqrt{1-k^2\sin^2y_2},\label{sdn}\\
\cos cx_1\,\sqrt{1-k^2\sin^2y_1}=\cos cx_2\,\sqrt{1-k^2\sin^2y_2}.\label{cdn}
\end{gather}

For generic $x_1,$ $y_1$ one has 
$\sin y_1\ne0$ and $\sin y_2\ne0.$
Hence we can divide equation~\eqref{ssn} by equation~\eqref{csn}
and obtain
$$
\tan ax_1=\tan ax_2\quad\Longleftrightarrow\quad x_1-x_2=\frac{\pi}{a}k,\quad k\in\mathbb{Z}.
$$
In the same way we obtain
$$
x_1-x_2=\frac{\pi}{b}l=\frac{\pi}{c}n,\quad l,n\in\mathbb{Z}
$$
from equations~\eqref{scn}--\eqref{cdn}. Hence we are looking for
integer $k,$ $l$ and $n$ such that
$$
\frac{k}{a}=\frac{l}{b}=\frac{n}{c}.
$$
Let $\frac{k}{a}=\frac{p}{q},$ where $(p,q)=1,$ $q>0.$
Then $k=\frac{pa}{q}$ and it follows
that $q$ divides $a.$ In the same way we prove that $q$ divides $b$ and $c$.
But we assume $(a,b,c)=1,$ hence $q=1$ and $k=pa,$ $l=pb,$ $n=pc.$ It follows
that
$$
x_1-x_2=\frac{\pi}{a}pa=\frac{\pi}{b}pb=\frac{\pi}{c}pc=p\pi,\quad p\in\mathbb{Z}.
$$
Since all functions in system~\eqref{ssn}--\eqref{cdn} are $2\pi$-periodic,
it is sufficient to consider the case $p=0$ and the case $p=1.$

In the case $p=0$ we have $x_2=x_1$ and system~\eqref{ssn}--\eqref{cdn}
implies the following system of equations,
$$
\sin y_1=\sin y_2,\quad \cos y_1=\cos y_2, \quad \sqrt{1-k^2\sin^2y_1}=\sqrt{1-k^2\sin^2y_2}.
$$
This system implies $y_2=y_1$ and we have only the trivial transformation $\Phi=id.$

In the case $p=1$ one has $x_2=x_1+\pi$ and
system~\eqref{ssn}--\eqref{cdn}
implies the following system of equations,
\begin{gather}
\sin y_1=(-1)^a\sin y_2,\label{sneq}\\
\cos y_1=(-1)^b\cos y_2,\label{cneq}\\
\sqrt{1-k^2\sin^2y_1}=(-1)^c\sqrt{1-k^2\sin^2y_2}.\label{dneq}
\end{gather}

If $a$ is even then equation~\eqref{sneq} implies that
either $y_2=y_1$ or $y_2=\pi-y_1.$ But $y_2=y_1$
implies from equations~\eqref{cneq} and~\eqref{dneq}
that $b$ and $c$ are even. However, this contradicts
our assumption $(a,b,c)=1.$ If $y_2=\pi-y_1$
then by equations~\eqref{cneq} and~\eqref{dneq}
we obtain that $b$ is odd and $c$ is even. Direct
check shows that in fact if $a$ and $c$ is even and $b$
is odd then the transformation $\Phi_1(x,y)=(x+\pi,\pi-y)$
satisfies equation~\eqref{Phi}.

If $a$ is odd then equation~\eqref{sneq} implies that
either $y_2=-y_1$ or $y_2=y_1+\pi.$ If $y_2=-y_1$
then equations~\eqref{cneq} and~\eqref{dneq}
imply that $b$ and $c$ are even. Direct
check shows that in fact if $a$ is odd and $b$ and
$c$ are even then the transformation $\Phi_2(x,y)=(x+\pi,-y)$
satisfies equation~\eqref{Phi}. If $y_2=y_1+\pi$
then equations~\eqref{cneq} and~\eqref{dneq}
imply that $b$ is odd and $c$ are even.
Direct check shows that in fact if $a$ and $b$
are odd and $c$ is even then the transformation $\Phi_3(x,y)=(x+\pi,y+\pi)$
satisfies equation~\eqref{Phi}.

These transformations imply the statement 5) of Theorem~\ref{maintheorem}.
The transformations $\Phi_1$ and $\Phi_2$ coincide under isometry
$T_{a,b,c}\cong T_{b,a,c}$ and correspond to the case~I). Due to
the isometry $T_{a,b,c}\cong T_{b,a,c}$ we can consider only
the case of odd $a.$ Then the points of $T_{a,b,c}$ have
unique coordinates $0\leqslant x<\pi$ and $-\pi\leqslant y<\pi$
and functions on $T_{a,b,c}$ could be considered as two-periodic
functions on $\mathbb{R}^2$ of period $2\pi$ with additional
invariance with respect to the transformation $\Phi_2.$

The transformation $\Phi_3$ corresponds to the case II).
The functions on $T_{a,b,c}$ could be considered as two-periodic
functions on $\mathbb{R}^2$ of period $2\pi$ with additional
invariance with respect to the transformation $\Phi_3.$

The case where there is no transformation $\Phi$ corresponds to the
case III). The functions on $T_{a,b,c}$ could be considered 
as two-periodic functions on $\mathbb{R}^2$ of period~$2\pi.$

The induced metric $g$ on $T_{a,b,c}$ is given by the
formula
\begin{equation}\label{g_T}
g=P(y)\,dx^2+%
\frac{2P(y)}{Q+2P(y)}\,dy^2,
\end{equation}
where
$$
P(y)=\frac{1}{2}(c^2+(b^2-a^2)\cos2y),\quad Q=c^2-a^2-b^2.
$$
The area of $T_{a,b,c}$ is obtained by a straightforward calculation.
Remark that $S(a,b,c)=S(b,a,c)$ since $T_{a,b,c}\cong T_{b,a,c}.$

Next, we should find $N(2)$ in order to investigate
extremal spectral properties of surfaces $T_{a,b,c}$
using El Soufi and Ilias Theorem~\ref{ElSoufi-Ilias}.
We follow the approach proposed in the paper~\cite{Penskoi2012}
and futher developed in the papers~\cite{Penskoi2013} 
and~\cite{Karpukhin2013,Karpukhin2013q1}. In order to shorten
the text we omit some
details and refer the reader to the paper~\cite{Penskoi2012}.

As it was explained before, the surfaces $T_{a,b,c}$
were constructed in such a way that one has a separation
of variables in the spectral problem for the
Laplace-Beltrami operator $\Delta\psi=\lambda\psi.$
More precisely, since $\Delta$ commutes with
$\frac{\partial}{\partial x},$ one can look for an eigenfunction
basis consisting of functions of the form
\begin{equation}\label{psi-form}
\psi(x,y)=\varphi(y)\sin lx\quad\mbox{or}\quad 
\psi(x,y)=\varphi(y)\cos lx,
\end{equation}
where we consider $l$ as an integer parameter in $\varphi(y).$
Substituting functions~\eqref{psi-form} in $\Delta\psi=\lambda\psi$
and separating variables one obtains the following equation,
% $$
% -\frac{\sqrt{Q+2P(y)}}{2P(y)}\left(\sqrt{Q+2P(y)}\varphi'(y)\right)'+%
% \left(\frac{l^2}{P(y)}-\lambda\right)\varphi(y)=0.
% $$
\begin{equation}\label{SL-equation}
\left(1+\frac{Q}{2P(y)}\right)\varphi''(y)+\frac{P'(y)}{2P(y)}\varphi'(y)+%
\left(\lambda-\frac{l^2}{P(y)}\right)\varphi(y)=0.
\end{equation}
The case of Lawson surfaces $\tau_{m,n}$
corresponds to $Q=0$ and was studied in the
paper~\cite{Penskoi2012}. We denote occasionally a solution
of~\eqref{SL-equation} by $\varphi(y,l)$ when
we need to emphasize dependence on the parameter $l.$

In the subcase III) the conditions
of $2\pi$-periodicity of $\psi(x,y)$ impose the $2\pi$-periodicity
condition on $\varphi(y).$ Thus, we should consider the periodic
Sturm-Liouville problem consisting of equation~\eqref{SL-equation}
and the periodicity condition
\begin{equation}\label{SL-periodicity}
\varphi(y+2\pi)\equiv\varphi(y).
\end{equation}
For each $l=0,1,2,3\dots$ we obtain the spectrum 
$\lambda_0(l)<\lambda_1(l)\leqslant\lambda_2(l)<\dots$
of the periodic Sturm-Liouville problem~\eqref{SL-equation},
\eqref{SL-periodicity}. Each eigenvalue $\lambda_i(l)$
corresponds to two eigenvalues $\lambda_k=\lambda_{k+1}$
of the initial Laplace-Beltrami
operator $\Delta$ on $T_{a,b,c}$ since one eigenfunction
$\varphi(y)$ of the problem~\eqref{SL-equation},
\eqref{SL-periodicity} correponds to two eigenfunctions~\eqref{psi-form}
of $\Delta.$ The only exception is the case of $l=0.$
Since $\sin 0x\equiv0,$ one eigenvalue
$\lambda_i(0)$ correponds to exactly one eigenvalue $\lambda_k$
of $\Delta.$
This implies that
\begin{equation}\label{N-formula}
N(2)=\#\{\lambda_i<2\}=\#\{\lambda_i(0)<2\}+2\#\{\lambda_i(l)<2,l\geqslant1\}.
\end{equation}
According to the Sturm oscillation theorem, one has the inequality
\begin{equation}\label{rows}
\lambda_0(l)<\lambda_1(l)\leqslant\lambda_2(l)<\lambda_3(l)\leqslant\lambda_4(l)<\dots
\end{equation}
On the other hand, one has the inequality
\begin{equation}\label{columns}
\lambda_i(0)<\lambda_i(1)<\lambda_i(2)<\lambda_i(3)<\lambda_i(4)<\dots.
\end{equation}
A simple proof of this inequality could be found in the paper~\cite{Karpukhin2013q1}.
The initial argument in the paper~\cite{Penskoi2012} uses the particular
properties of Lawson tau-surfaces and quite complicated.

Let us now remark that we know three eigenvalues of the problem~\eqref{SL-equation},
\eqref{SL-periodicity} equal to $2.$ Indeed, Takahashi Theorem~\ref{takahashi}
states that the components of the immersion~\eqref{Fabc}
are eigenfunctions of $\Delta$ with eigenvalue $2.$ They are exactly of 
the form~\eqref{psi-form} with $l=a,$ $l=b$ and $l=c.$
Let us look at $\sin cx\,\tilde\varphi_3(y)$ and $\cos cx\,\tilde\varphi_3(y).$
Since they are eigenfunctions of $\Delta,$ we know that
$\tilde\varphi_3(y)$ is an eigenfunction with eigenvalue $2$
of the problem~\eqref{SL-equation},
\eqref{SL-periodicity} with $l=c.$ Let us remark that $\tilde\varphi_3(y)$
has no zeroes. Then Sturm oscillation theorem implies that this is
an eigenfunction corresponding to $\lambda_0(c).$ Hence, $\lambda_0(c)=2.$
In a similar way we can establish that if $a\geqslant b$ then $\lambda_1(a)=2$
and $\lambda_2(b)=2$ and if $a<b$ then $\lambda_1(b)=2$
and $\lambda_2(a)=2,$ see the paper~\cite{Penskoi2012} for more details.

It follows now from inequality~\eqref{columns} that among
all eigenvalues $\lambda_0(l)$ exactly $\lambda_0(0),\dots,\lambda_0(c-1)$
are less than $2.$ The similar statement holds for $\lambda_1(l)$
and $\lambda_2(l).$

Using the theory of the Lam\'e equation one can prove that
$\lambda_3(0)>2,$ see the papers~\cite{Penskoi2012} and~\cite{Karpukhin2013q1}
for more details. Then inequality~\eqref{columns} implies
that for any $l$ we have $\lambda_3(l)>2.$ Then inequality~\eqref{rows} implies
that for any $i>3$ and any $l$ we have $\lambda_i(l)>2.$ Hence,
we can find $N(2)$ by formula~\eqref{N-formula} and obtain
$$
N(2)=3+2(a-1+b-1+c-1)=2(a+b+c)-3.
$$
This means that in the subcase III) the metric~\eqref{g_T}
induced on the torus $T_{a,b,c}$ is extremal for
the functional $\Lambda_j(\mathbb{T}^2,g),$
where $j=2(a+b+c)-3.$ The corresponding value
of this functional $\Lambda_j(T_{a,b,c})=2\Area(T_{a,b,c})=2S(a,b,c).$
This proves the part of the statement 6) of theorem~\ref{maintheorem}
concerning the subcase III). The special cases of $T_{a,0,c}$ and
$T_{0,0,1}$ could be investigated in the same way.

The subcases I) and II) are more complicated since one have to take
into account not only the $2\pi$-periodicity conditions~\eqref{SL-periodicity}.
In the subcase I) eigenfunctions $\psi(x,y)$ have to satisfy also
the condition of invariance with respect to the transformation~$\Phi_2.$
Since we look for eigenfunctions of the form~\eqref{psi-form}, this
condition could be written in the following way: $\varphi(x,l)$
has to be even function for even $l$ and odd function
for odd $l.$ One can then find $N(2)$ in the same way as
before but taking into the account the parity of solutions.
All details could be found in the proof of the 
Main Theorem in the paper~\cite{Penskoi2012} and we give
here only the answer,
$$
N(2)=a+b+c-3.
$$
This means that in the subcase I) the metric~\eqref{g_T}
induced on the Klein bottle $T_{a,b,c}$ is extremal for
the functional $\Lambda_j(\mathbb{KL},g),$
where $j=a+b+c-3.$ The corresponding value
of this functional $\Lambda_j(T_{a,b,c})=2\Area(T_{a,b,c})=S(a,b,c).$
This proves the part of the statement 6) of theorem~\ref{maintheorem}
concerning the subcase I). The special case of $T_{a,0,c}$
could be investigated in the same way.

In the subcase II) eigenfunctions $\psi(x,y)$ have to satisfy also
the condition of invariance with respect to the transformation~$\Phi_3.$
This condition means that $\varphi(y,l)$ has to be $\pi$-periodic
for even $l$ and $\pi$-antiperiodic for odd $l.$
One can then find $N(2)$ in the same way as
before but taking into the account the $\pi$-(anti)periodicity of solutions.
All details could be found in the proof of the 
Main Theorem in the paper~\cite{Penskoi2012} and we give
here only the answer,
$$
N(2)=a+b+c-3.
$$
This means that in the subcase II) the metric~\eqref{g_T}
induced on the torus $T_{a,b,c}$ is extremal for
the functional $\Lambda_j(\mathbb{T}^2,g),$
where $j=a+b+c-3.$ The corresponding value
of this functional $\Lambda_j(T_{a,b,c})=2\Area(T_{a,b,c})=S(a,b,c).$
This proves the part of the statement 6) of theorem~\ref{maintheorem}
concerning the subcase II).

This finishes the proof. $\Box$

We should remark that after the author's talk at the Analisys
Seminar at the McGill University I.~Polterovich conjectured
that one can find tori minimally immersed in spheres
using not only three first solutions $\dn z,$ $\cn z$
and $\sn z$ of the Lam\'e equaiton with $n=1,$ but also
using next solutions. Is not clear how
to prove this conjecture since next solutions are
given only by series.

\section{The Klein bottle $T_{1,0,2}$}\label{T102-section}

It follows from Theorem~\ref{maintheorem}
that $T_{1,0,2}$
is a Klein bottle and the metric on $T_{1,0,2}$
is extremal for $\Lambda_1(\mathbb{KL},g).$
In the same time, El~Soufi, Giacomini and Jazar proved
in paper~\cite{ElSoufi-Giacomini-Jazar2006} 
that the metric on $\tilde{\tau}_{3,1}$ is the unique
(up to multiplication by a constant)
extremal metric for the first eigenvalue
on the Klein bottle and hence the maximal one.

It is interesting to compare the values $\Lambda_1(T_{1,0,2})$
and $\Lambda_1(\tilde{\tau}_{3,1}).$ We have 
$$
\Lambda_1(T_{1,0,2})=S(1,0,2)=S(0,1,2)=%
2\pi\left(8E\left(\frac{1}{2}\right)%
-3K\left(\frac{1}{2}\right)\right),
$$
where $S(a,b,c)=S(b,a,c)$ since $T_{a,b,c}\cong T_{b,a,c},$
and
$$
\Lambda_1(\tilde{\tau}_{3,1})=12\pi E\left(\frac{2\sqrt{2}}{3}\right).
$$
Both values are equal due to the identity
(see the book~\cite{BatemanErdelyi1955})
$$
E\left(\frac{2\sqrt{k}}{1+k}\right)=%
\frac{2E(k)-(k')^2K(k)}{1+k},
$$
where $k'=\sqrt{1-k^2}.$

Since in both cases $\lambda_1=2,$ the areas are equal. This implies
the following Proposition.

\begin{Proposition} The metric on the Klein bottle $T_{1,0,2}$
is maximal for $\Lambda_1(\mathbb{KL},g).$
The Klein bottle $T_{1,0,2}$ is isometric to the bipolar
Lawson Klein bottle $\tilde{\tau}_{3,1}.$
\end{Proposition}

It would be interesting to find this isometry explicitly.
The explicit parametrisation of $T_{1,0,2}\subset\mathbb{S}^4$ is given by the
formula
$$
\left(\frac{1}{\sqrt{2}}\sin x\sin y,\frac{1}{\sqrt{2}}\cos x\sin y,%
\sqrt{\frac{5}{8}}\cos y,\right.
$$
$$
\left.\sqrt{\frac{3}{8}}\sin2x\sqrt{1+\frac{1}{3}\sin^2y},%
\sqrt{\frac{3}{8}}\cos2x\sqrt{1+\frac{1}{3}\sin^2y}\right),
$$
where we omit one of the components equal to zero.
It would be also interesting to find whether
other bipolar Lawson surfaces $\tilde{\tau}_{m,k}$
are among $T_{a,b,c}.$

\section*{Acknowledgments}

The author thanks D.~Jakobson, I.~Polterovich
and P.~Winternitz for fruitful discussions
at the Centre de
Recherches Math\'ematiques, Univerisit\'e
de Montr\'eal (CRM).
The author is very grateful to the CRM for its hospitality.

The author also thanks M.~Karpukhin for useful
discussions.

This work was partially supported
by Russian Federation Government grant no.~2010-220-01-077, 
ag. no.~11.G34.31.0005,
by the Russian Foundation
for Basic Research grant no.~11-01-12067-ofi-m-2011,
by the Russian State Programme for the Support of
Leading Scientific Schools grant no.~4995.2012.1
and by the Simons-IUM fellowship.

\end{document}